\documentclass[12pt, a4paper]{article}

\usepackage[utf8]{inputenc}
\usepackage[T1]{fontenc}
\usepackage[english]{babel} 
\usepackage[margin=2.5cm]{geometry}

\usepackage{setspace}
\usepackage{lineno} 

\usepackage{amsmath, amssymb, amsthm}

\usepackage{graphicx}
\usepackage{booktabs} 
\usepackage{caption}
\usepackage{subcaption}

\usepackage[
    backend=biber,
    style=authoryear,    
    sorting=nyt,         
    maxcitenames=2,      
    mincitenames=1,      
    maxbibnames=99,      
    giveninits=true,     
    uniquename=false,    
    natbib=true          
]{biblatex}
\usepackage{authblk} 

\usepackage[colorlinks=true, allcolors=blue]{hyperref}

\title{Efficient Path Reconstruction in Prehistoric Human Migration: An Adaptive Dijkstra's Algorithm Based on Wavelet Compression for Topographic Data}

\author[1]{Max Brockmann, OrcidID: 0009-0007-1197-0634}
\author[1]{Lena Perlberg, OrcidID: 0009-0000-9494-0236}
\author[1]{Angela Kunoth, OrcidID: 0000-0001-5254-1501} 
\affil[1]{Department of Mathematics and Computer Science, University of Cologne, Germany}
\date{June 30, 2026}

\begin{document}

\maketitle

\begin{abstract}
The reconstruction of prehistoric migration routes requires the accurate evaluation of effective distances across complex topographies. Traditional Least-Cost Path Analyses (LCPA) using high-resolution elevation models often lead to dense computational grids, making algorithms like Dijkstra's prohibitively slow for large-scale archaeological modelling. In this paper, we present a novel approach using adaptive wavelet methods to dynamically compress topographic data. By retaining high resolution only in areas of high topographic complexity (e.g., mountain passes) and smoothing homogeneous regions, we significantly accelerate the Dijkstra algorithm without losing essential routing topology. We demonstrate the efficacy of this method using migration scenarios.
\end{abstract}

\section{Introduction}
Reconstructing the dynamics of past human migration and interaction has become a central objective in the interdisciplinary research in HESCOR. Because the physical movement of prehistoric populations rarely leaves direct archaeological traces, we rely on computational models to infer ancient mobility networks \citep{hilpert_interdisciplinary_mobility}. A foundational method for this task is a Least-Cost Path Analysis (LCPA), which calculates ``effective distances'' by accounting for environmental constraints such as mountain ranges, steep inclines, or water bodies \citep{herzog2014least}. 
\par
To achieve realistic reconstructions, archaeologists commonly evaluate high-resolution digital elevation models (DEMs) \citep{Conolly_Lake_2006}, such as the global Earth TOPOgraphy (ETOPO) dataset as shown in Figure \ref{fig:adaptive_mesh_example_a} \citep{noaa_etopo_2022}. These landscapes are translated into dense grid-based graphs, where the travel effort between nodes is estimated using empirical cost functions like Tobler's Hiking function. However, as the spatial resolution of topographic data increases, we face a severe computational bottleneck. Solving these routing problems across vast spatial datasets fundamentally relies on algorithms like Dijkstra's algorithm. While mathematically exact, the asymptotic time complexity of Dijkstra's algorithm scales non-linearly \citep{cormen2022introduction}. Consequently, evaluating continental-scale networks at high resolutions in general exceeds practical computational and memory capacities. 
\par
A natural workaround is uniform data compression, which arbitrarily reduces the grid resolution to accelerate processing times \citep{tang2023effective}. From a methodological standpoint, however, this uniform simplification structurally distorts the landscape. Crucial fine-scale topographic features that would have dictated human movement, such as narrow mountain passes or steep valley corridors, are frequently eradicated by uniform smoothing, leading to inherently flawed path reconstructions.
\par
Motivated by the computational challenges encountered by archaeological models within the interdisciplinary HESCOR project, this paper introduces a mathematical framework to overcome this topographic bottleneck. We propose an adaptive multi-scale routing approach based on the Fast Wavelet Transform. Rather than employing a static grid, our algorithm dynamically takes into account the local regularity of the terrain. By thresholding wavelet coefficients, the framework automatically merges topologically homogeneous regions into large, computationally inexpensive spatial blocks, while preserving maximum geometric resolution exclusively in rugged areas with high topographic variance. 
\par
We demonstrate how Dijkstra's algorithm can be structurally adapted to navigate this dynamic, multi-scale mesh directly. This paper focuses on the computational foundations and the algorithmic validation of this method, demonstrating how it massively reduces the memory footprint and runtime without compromising the topological integrity of the resulting pathways. Specifically, by leveraging the mathematical properties of second order B-spline Wavelets to achieve maximum spatial compression without sacrificing topological fidelity, our framework enables the computation of large-scale data that was previously unexplored. The efficient implementation of the adaptive Dijkstra's algorithm and its components are discussed in detail as part of the \texttt{ArcheoGra.jl} Julia-package \citep{brockmann_multigrid_2026}. A thorough investigation of the mathematical, theoretical properties of the adaptive Dijkstra's algorithm is provided in \citep{brockmann_2026_adaptive}.
\par
The results shown here are based on prior work on Wavelet compression and Least-Cost Paths for topographic data from \citep{perlberg2025wavelet} and the extension of Wavelet compression to adaptive piecewise linear representations with the corresponding adaptive Dijkstra's algorithm from \citep{brockmann_multigrid_2026}.

\section{Topography and Effective Distance in Archaeological Modelling}
\label{sec:topography_and_lcpa}

In spatial mobility models, evaluating the geographic relationship between distant locations using simple geometric measures, such as the Haversine great-circle distance, is often a strong simplification \citep{herzog2014least}. Geometric distances fail to account for environmental and topographic constraints that drastically influence human movement. Geographic hurdles, such as mountain ranges, steep valleys, or impassable water bodies, can significantly increase the actual travel effort between two locations, rendering a geometrically short route practically unviable. 
\par
To compute realistic migration pathways, the geometric distance can be replaced with an ``effective distance'' by framing the routing task as a Least-Cost Path Analysis (LCPA) problem \citep{herzog2014least}. LCPA translates the physical landscape into a mathematical graph structure, allowing routing algorithms to find the optimal trajectory that minimizes the total travel cost, where cost can be energy expenditure or travel time.

\subsection{Discretizing the Topographic Landscape}
To apply LCPA computationally, the continuous landscape must first be discretized into a finite routing graph, denoted as $\mathcal{G} = (\mathcal{V}, \mathcal{E})$. Here, the vertex set $\mathcal{V}$ represents distinct spatial locations, and the edge set $\mathcal{E}$ defines the possible movement transitions between them. 
\par
For global and continental-scale modelling, we utilize the Earth TOPOgraphy \\(ETOPO) dataset, which provides a high-resolution relief model of the Earth's surface \citep{noaa_etopo_2022}. Specifically, we operate on a 60 arc-second resolution grid, which equates to spatial cells of approximately $1 \times 1$ km to $1.5 \times 1.5$ km. In our graph representation, the geometric centre of each grid cell serves as a vertex $v \in \mathcal{V}$, carrying the local elevation value. To restrict human movement to realistic landmasses, we apply geospatial water masks to identify and remove all vertices corresponding to large oceanic bodies or lakes from the routing domain.
\par
To enable movement across this terrain, we construct an 8-connected grid graph. Each vertex is connected via edges to its eight immediate horizontal, vertical, and diagonal neighbours. This topology allows the routing algorithm to navigate the landscape with a high degree of spatial flexibility, while keeping the ``one-step'' reach of the allowed connections realistic. 

\subsection{Translating Slope into Travel Cost}
In order for the routing algorithm to evaluate different pathways, each edge in the graph must be assigned a mathematical weight or ``cost'' that reflects the physical difficulty of traversing that specific segment. Here, we consider topographic slope $S$ as the primary determinant of travel effort over large distances.
\par
To quantify this effort, we employ empirical models based on Tobler's Hiking function. In this framework, we utilize the modified hiking function by \citet{kondo2010gps}. This function realistically captures the asymmetry of human movement, where walking speed $KS(S)$ in $\text{km/h}$ depends on the slope $S$:
\begin{equation}
    KS(S) := \begin{cases}
        5.1 \cdot \exp\left(-2.25 \cdot |S + 0.07|\right), & \text{if } S \geq -0.07 \text{ (uphill or slight downhill)}\\
        5.1 \cdot \exp\left(-1.5 \cdot |S + 0.07|\right), & \text{if } S < -0.07 \text{ (steep downhill)},
    \end{cases}
\end{equation}
where the different cases account for the directional asymmetry of the effort.
\par
The fastest walking speed is achieved at a slight downhill slope of $-7\%$. For the algorithmic implementation, the traversal cost $c(u,v)$ assigned to a directed edge $(u,v)$ is defined as the travel time, which is the physical distance divided by the walking speed $KS(S)$. Note that the cost function is asymmetric ($c(u,v) \neq c(v,u)$), which impacts the shortest-path computation.

\subsection{The Computational Bottleneck of LCPA}
Once the landscape is represented as a weighted graph, finding the path of minimum total cost is achieved using a Least-Cost Path algorithm. Least-Cost Paths evaluate the cumulative cost of traversing from a starting vertex to a target vertex, iteratively exploring the graph to identify the optimal route. A representation of a (non-adaptive) LCPA is shown in Figure \ref{fig:alpine_paths_comparison} in the black paths. While there are several possible LCPA algorithms (notably $A^*$), see for an example \citep{Bast2016}, for the application within the HESCOR project, we require an ``All-Pairs Shortest Paths'' (APSP) approach. For this specific application, Dijkstra's algorithm \citep{dijkstra_note_1959} is chosen. It guarantees the mathematical discovery of the exact shortest path by iteratively evaluating the cost of reaching all accessible vertices from a given starting point. For a detailed discussion of the algorithmic implementation and its resulting limitations see \citep{brockmann_2026_adaptive}.
\par
However, applying this standard algorithm to high-resolution topographic data reveals a severe computational bottleneck. The algorithmic efficiency of Dijkstra's algorithm, as measured using Big-$\mathcal{O}$ notation, describes how the runtime scales with the size of the input. Even when implemented efficiently, Dijkstra's algorithm exhibits a time complexity of $\mathcal{O}(|\mathcal{E}| + |\mathcal{V}| \log |\mathcal{V}|)$, where $|\mathcal{V}|$ is the number of vertices and $|\mathcal{E}|$ is the number of edges \citep{cormen2022introduction}, meaning that any increase in the number of vertices increases non-linearly in the $|\mathcal{V}| \log |\mathcal{V}|$ term, as well as the $|\mathcal{E}|$ term, as any additional vertices will also increase the number of edges.
\par
In practical terms, this relationship creates a severe computational bottleneck when working with high-resolution topographic data. If one doubles the spatial resolution of a digital elevation model (e.g., from 1 km to 0.5 km per grid cell), the number of vertices $|\mathcal{V}|$ in the 2D grid quadruples. Because the algorithm’s runtime grows faster than linearly with the number of nodes, this refinement leads to a non-linear increase in computation time that significantly exceeds the initial quadrupling of data points. When simulating continental-scale networks between hundreds of archaeological sites, this scaling behaviour renders standard LCPA on dense, uniform grids computationally prohibitive.

\section{Methodology: Adaptive Wavelet Compression}
\label{sec:adaptive_wavelet_compression}

To resolve the topographic bottleneck inherent in standard LCPA, we abandon the constraint of a uniform grid. The core objective is to create a dynamic spatial representation that allocates high resolution only where it is strictly necessary, such as in rugged mountains, while aggressively compressing topographically monotonous regions like vast plains. To achieve this mathematically, we employ the Fast Wavelet Transform (FWT) \citep{dahmen_wavelet_1997, dahmen1999biorthogonal}. Wavelet methods and adaptive Dijkstra's algorithms have previously been applied to topographic data to achieve a fast framework for LCP in \citep{cowlangi2008multiresolution,tsiotras2012multiresolution}, however not yet in archaeological application. Furthermore, current studies limited the scope of Wavelet compression to a piecewise constant representation of the underlying data. Here, we will extend this scheme to piecewise linear representations of data, to further improve the compression efficacy and thus allowing for an increased speed-up of the adaptive Dijkstra's algorithm. A complete discussion of the theoretical properties of the FWT in application to topographic data and Dijkstra's algorithm is provided in \citep{brockmann_2026_adaptive}.
\par
\begin{figure}[htbp]
    \centering
    \begin{subfigure}{0.45\textwidth}
        \centering
        \includegraphics[width=0.8\textwidth]{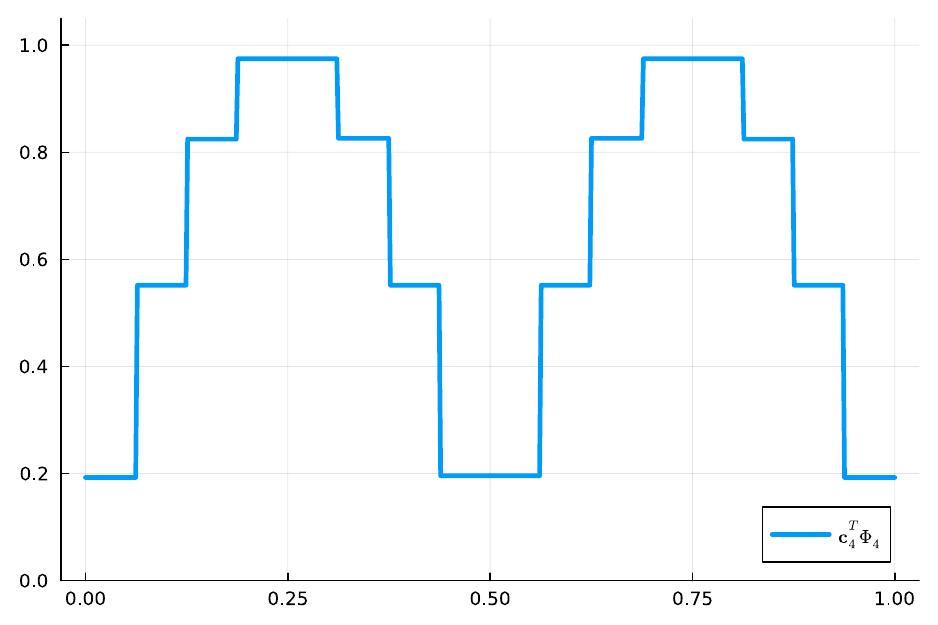}
        \caption{Piecewise constant approximation $N_1$ of the trial function $f(x):=|\sin(2\pi x)|$ using 16 basis functions.}
        \label{fig:piecewise_approximations_a}
    \end{subfigure}
    \hfill
    \begin{subfigure}{0.45\textwidth}
        \centering
        \includegraphics[width=0.8\textwidth]{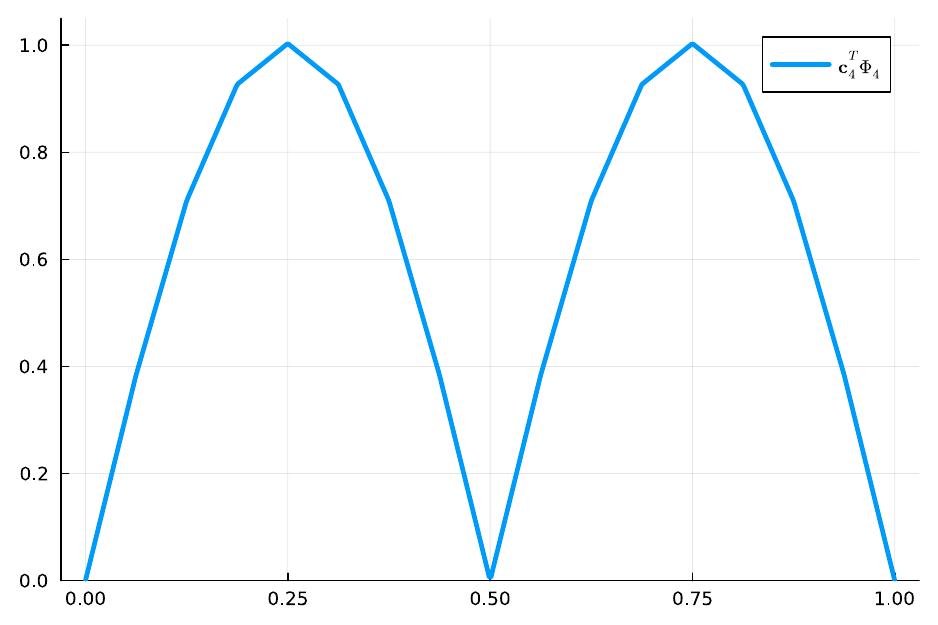}
        \caption{Piecewise linear approximation $N_2$ of the trial function $f(x):=|\sin(2\pi x)|$ using 16 basis functions.}
        \label{fig:piecewise_approximations_b}
    \end{subfigure}
    \caption{Comparison of piecewise constant ($N_1$) and piecewise linear ($N_2$) approximations of the trial function $f(x):=|\sin(2\pi x)|$. The $N_2$ approximation provides a smoother and more accurate representation of the underlying function. (Figure from \citep{perlberg2025wavelet})}
    \label{fig:piecewise_approximations}
\end{figure}

We refer to piecewise constant functions as $N_1$, and continuous piecewise linear functions as $N_2$. The properties of the basis functions are reflected in their approximation of underlying data. In Figure \ref{fig:piecewise_approximations} we show the approximation of the (one-dimensional) trial function $f(x):=|\sin(2\pi x)|$ using $N_1$ in Figure \ref{fig:piecewise_approximations_a} and $N_2$ in Figure \ref{fig:piecewise_approximations_b}. The piecewise constant approximation $N_1$ creates a blocky representation of the underlying function, while the piecewise linear approximation $N_2$ creates a continuous representation of the underlying function. The increased smoothness (piecewise linear vs.\ piecewise constant) of the $N_2$ basis functions allows for a more accurate representation of the underlying data.
\par
The core mechanism of wavelet theory is the multi-scale decomposition of a function, in our case, the topographic elevation function $f(x, y)$. Instead of merely storing the raw elevation of each pixel, which we refer to as the single-scale representation of data defined as 
\[
    f_L(x,y) = \sum_{k} c_{L,k} \phi_{L,k}(x,y),
\]
the discrete wavelet transform represents the landscape as a projection onto nested vector spaces $\mathcal{S}_0 \subset \mathcal{S}_1 \subset \dots \subset \mathcal{S}_L$, where $L$ is the finest resolution level. The topographic function $f_L$ at maximum resolution can be mathematically expressed via the multi-scale wavelet decomposition:
\begin{equation}\label{eq:wavelet_decomposition}
f_L(x,y) = \sum_{k} c_{0,k} \phi_{0,k}(x,y) + \sum_{\ell=0}^{L-1} \sum_{k} d_{\ell,k} \psi_{\ell,k}(x,y).
\end{equation}
Here, the scaling functions $\phi_{0,k}$ and their coefficients $c_{0,k}$ define the coarse baseline approximation of the landscape. The crucial elements are the wavelet basis functions $\psi_{\ell,k}$ and their corresponding detail coefficients $d_{\ell,k}$ \citep{dahmen_wavelet_1997}. These detail coefficients encode the geometric differences or topographic variance between adjacent spatial scales $\ell$ and $\ell+1$. 
For compressing data using the wavelet transform, the multi-scale representation allows us to identify and discard coefficients that fall below a certain threshold, effectively keeping only the ``most important'' topographic details. The exact mechanism for this compression is described in the following Section \ref{sec:local_regularity_and_hierarchical_validation}.
The advantage of wavelet compression is that the representation \eqref{eq:wavelet_decomposition} can be calculated in $\mathcal{O}(N)$, i.e. linear time \citep{dahmen_wavelet_1997}, where $N$ is the number of pixels in the original raster. This allows us to find the ideal setup for compression and detail retention in a computationally efficient manner, even for large-scale topographic datasets. Details of the highly efficient implementation of the Fast Wavelet Transform for 2D data can be found in \citep{perlberg2025wavelet}.
\par
\begin{figure}[htbp]
    \centering
    \begin{subfigure}[t]{0.45\textwidth}
        \centering
        \includegraphics[width=0.8\textwidth]{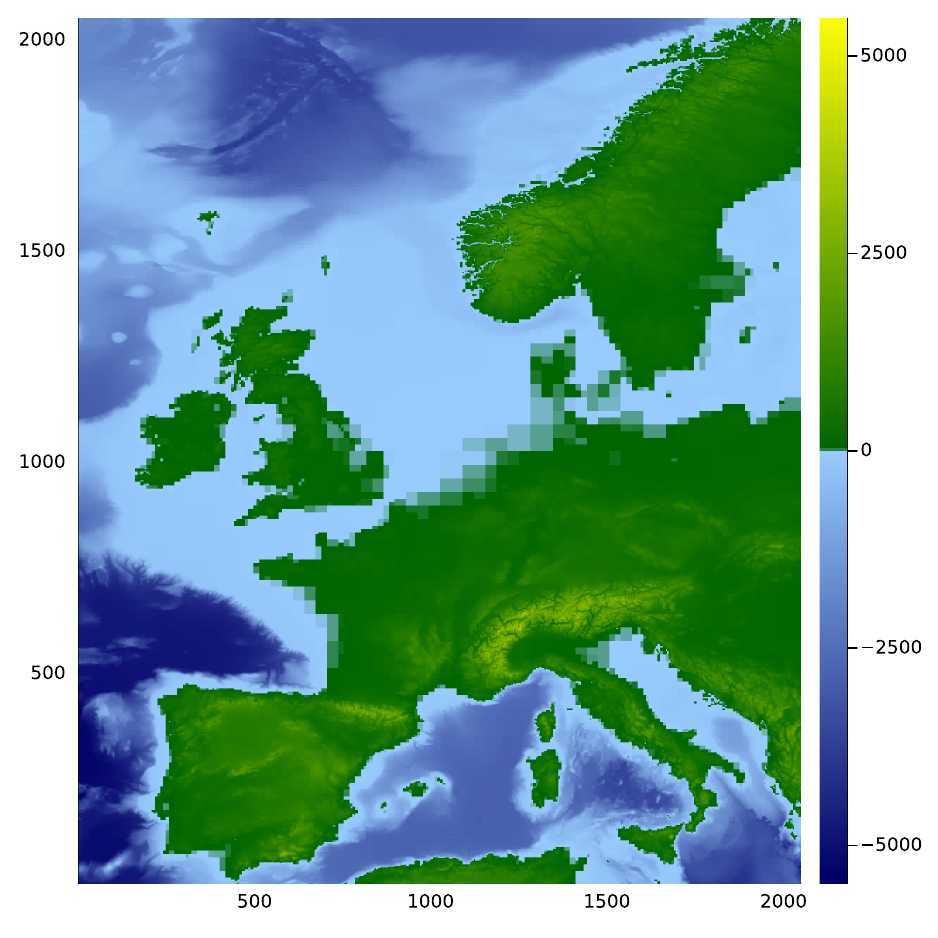}
        \caption{Adaptive Wavelet representation using $N_1$ basis functions and Haar wavelets of topographic data. The representation is compressed to $N=65536$ coefficients, resulting in a compression rate of $98.43\%$, meaning under $2\%$ of the original data is retained.}
        \label{fig:adaptive_mesh_example_a}
    \end{subfigure}
    \hfill
    \begin{subfigure}[t]{0.45\textwidth}
        \centering
        \includegraphics[width=0.8\textwidth]{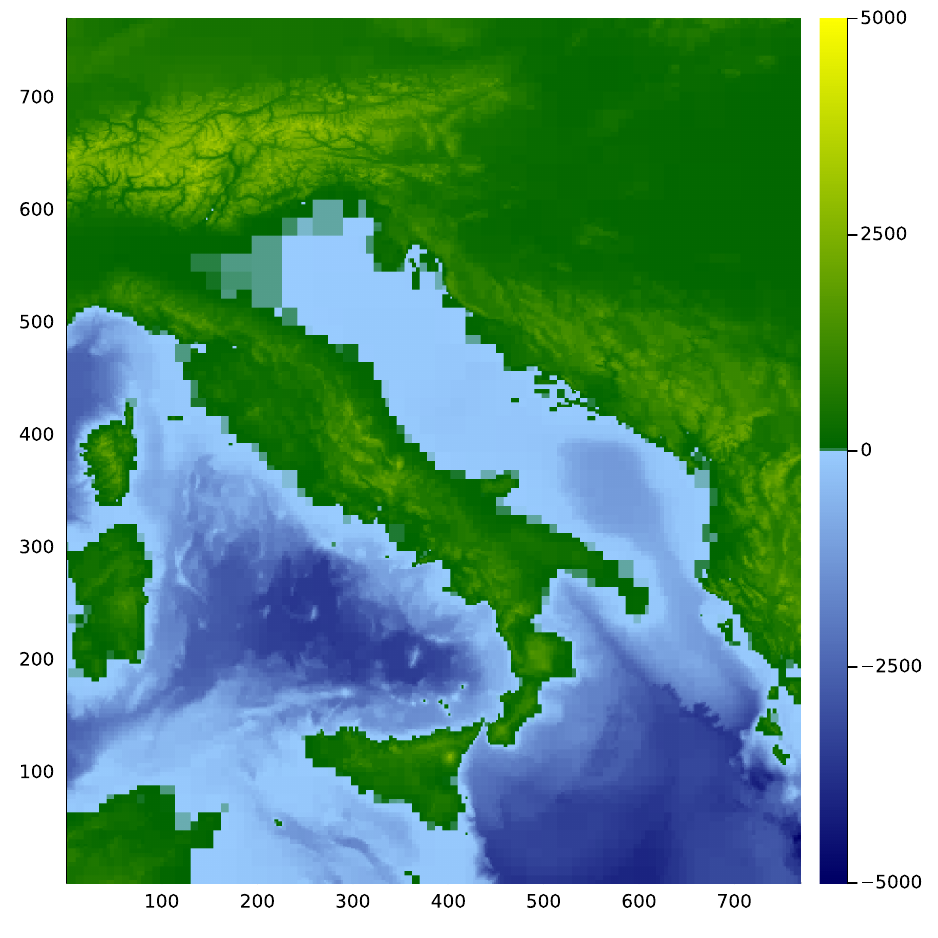}
        \caption{Cutout of the adaptive Wavelet representation from northern Italy, including the eastern parts of the Alps. The adaptive mesh preserves the narrow valleys and mountain passes, while compressing the flat coastal areas into large macroscopic blocks.}
        \label{fig:adaptive_mesh_example_b}
    \end{subfigure}
    \caption{Example of the adaptive representation of topographic data. The representation shows elevation of different regions as by colour intensity. In homogeneous regions, coefficients are disregarded and summarized in large spatial chunks. In areas of high fidelity, coefficients are retained and topography is represented in high resolution. (Figure adapted from \citep{perlberg2025wavelet})}
    \label{fig:adaptive_mesh_example}
\end{figure}

In the adaptive wavelet representation, the retained detail coefficients $d_{\ell,k}$ after compression, are the only information kept in the adaptive-resolution mesh. Any scaling or wavelet function $\phi_{\ell,k}$ or $\psi_{\ell,k}$ whose corresponding coefficient is set to zero during compression is discarded, and the corresponding spatial region is represented as a single macroscopic block. The resulting representation of data, consists of cells of varying sizes, as shown in Figure \ref{fig:adaptive_mesh_example}. Whereas in full representation every pixel is represented as a vertex, in the adaptive mesh, large homogeneous regions are represented by single vertices, drastically reducing the total number of vertices $|\mathcal{V}|$ that must be processed by Dijkstra's algorithm.

\subsection{Local Regularity and Hierarchical Validation}\label{sec:local_regularity_and_hierarchical_validation}
The primary structural advantage of this multi-scale formulation lies in evaluating the terrain's local regularity. In topographically homogeneous regions, the geometric variation between scales is minimal, causing the corresponding wavelet detail coefficients $d_{\ell,k}$ to drop near zero. To compress this data, we apply a Best-$N$-term thresholding strategy \citep{Cohen2003numerical}. This approach ranks the detail coefficients by magnitude and strictly retains only the $N$ largest ones, setting all remaining coefficients to zero. This targeted truncation safely discards redundant data points in flat areas while perfectly preserving the critical, large detail coefficients that represent steep slopes, ridges, or crucial mountain passes.
\par
However, translating this compression into a spatial routing mesh requires caution. A purely local evaluation might deem a large spatial block ``flat'' on average, prematurely consolidating it and thereby accidentally erasing a narrow, sub-scale gorge crucial for human movement. To prevent the distortion of sub-scale topographical barriers, we implement a bottom-up hierarchical validation scheme. In this framework, a macroscopic region is only classified as uniform, and thus eligible for consolidation, if all of its constituent finer sub-regions strictly lack significant topographic details.

\subsection{Adaptive Quadtree Partitioning and Basis Selection}
This hierarchical validation naturally translates the landscape into a quadtree data structure. Homogeneous regions are systematically grouped together, merging four smaller cells into one larger ``leaf cell'' at the next coarser level. By representing large flat areas as single nodes, the size of the vertex set $|\mathcal{V}|$ processed by Dijkstra's algorithm is drastically reduced.
\par
Translating this structure into a traversable graph requires choosing an appropriate mathematical basis function. Piecewise-constant functions ($N_1$ B-Splines / Haar wavelets) treat the landscape as mutually disjoint, flat building blocks. While computationally efficient, they create artificial, cliff-like steps at boundaries between different scales. Haar-wavelets have been successfully applied to topographic data, albeit non archaeologically, to achieve a fast framework for LCP in \citep{tsiotras2012multiresolution} and \citep{cowlangi2008multiresolution}. 
\par
To avoid these artifacts, we extend the framework to continuous piecewise-linear functions ($N_2$ B-Splines / Hat wavelets). Rather than isolated blocks, nodes in a Hat wavelet mesh represent the peaks of overlapping, tent-shaped spatial supports. This creates a significantly smoother, continuous representation of the terrain, yielding a more realistic translation of continuous topography for pathfinding algorithms.

\section{The Adaptive Dijkstra Algorithm}
\label{sec:the_adaptive_dijkstra_algorithm}

Transitioning from a uniform raster to a multi-scale wavelet mesh requires fundamental modifications to the routing mechanism. 
While we continue to employ the ``standard'' Dijkstra algorithm to compute the shortest path, the adaptive nature of the approach lies in the dynamic construction of the routing graph and its cost function. The adaptive graph consists of vertices representing vastly different spatial extents, from $1.5 \times 1.5$ km cells to blocks spanning tens of kilometers. Navigating this irregular topology, where a micro-scale valley cell might directly border a macro-scale plain, requires rigorous scale-aware edge evaluation.

\subsection{Graph Traversal Across Hierarchical Boundaries}
In a standard uniform grid, the edge set $\mathcal{E}$ is trivially defined by immediate spatial neighbours. In our adaptive multi-scale framework, however, a highly resolved leaf cell at maximum resolution $L$ may directly border a compressed macroscopic block at level $\ell < L$. Consequently, the routing algorithm must dynamically resolve transitions across these hierarchical boundaries.
\par
The definition of traversable edges strictly depends on the chosen multiresolution basis. For piecewise-constant Haar wavelets, adjacent cells of varying scales create ``hanging nodes.'' Transitioning across these boundaries requires an adaptive adjacency mapping that correctly connects multiple fine-scale vertices to a single coarse-scale vertex. The connectivity definition for piecewise constant wavelets is based on \citep{cowlangi2008multiresolution}.
\par
\begin{figure}[htbp]
    \centering
    \begin{subfigure}[t]{0.45\textwidth}
        \centering
        \includegraphics[width=0.6\textwidth]{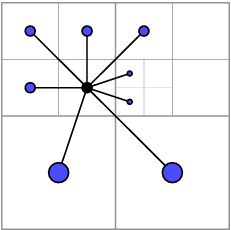}
        \caption{Edge definition for piecewise-constant Haar wavelets. An edge exists between two vertices if their corresponding basis functions share a common boundary. From the compressed wavelet representation, the connectivity structure finds connections between disjoint blocks of varying scales.}
        \label{fig:haar_edge_definition}
    \end{subfigure}
    \hfill
    \begin{subfigure}[t]{0.45\textwidth}
        \centering
        \includegraphics[width=0.6\textwidth]{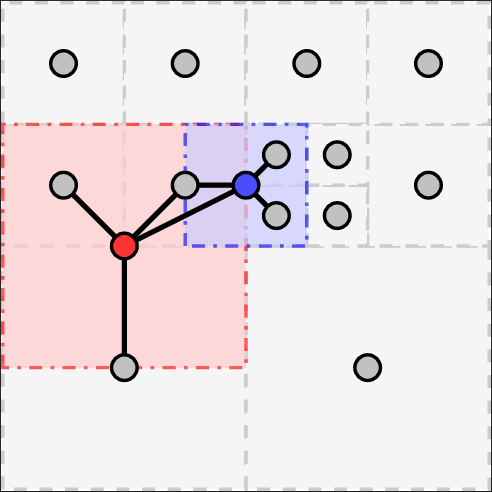}
        \caption{Edge definition for piecewise-linear $N_2$ B-splines. An edge exists between two vertices if their corresponding basis functions have overlapping support. The overlapping supports of the hat functions create a more interconnected graph structure compared to the disjoint blocks of Haar wavelets. The figure does not capture all supports of hat functions that would be contained in the full graph for visualization purposes, but only a subset to illustrate the concept.}
        \label{fig:hat_edge_definition}
    \end{subfigure}
    \caption{Edge definitions for different wavelet bases in the adaptive multi-scale routing framework. The connectivity rules ensure that the routing algorithm can navigate across hierarchical boundaries while preserving topological fidelity. (Figure from \citep{brockmann_multigrid_2026})}
    \label{fig:edge_definitions}
\end{figure}

Since for piecewise-constant Haar wavelets, the nodes represent disjoint blocks, an edge $(u,v)$ between two multi-scale vertices $u$ (at level $\ell_u$ with) and $v$ (at level $\ell_v$) is topologically defined if and only if their respective basis functions share a common boundary, as visualized in Figure \ref{fig:haar_edge_definition}. This approach has been developed in \citep{brockmann_multigrid_2026}.
\par
Conversely, for piecewise-linear $N_2$ B-splines, the nodes do not represent disjoint blocks but rather overlapping spatial supports. An edge $(u,v)$ between two multi-scale vertices $u$ (at level $\ell_u$ with translation $k_u$) and $v$ (at level $\ell_v$ with translation $k_v$) is topologically defined if and only if their respective basis functions share an overlapping support:
\begin{equation}
\operatorname{supp}(\psi_{\ell_u, k_u}) \cap \operatorname{supp}(\psi_{\ell_v, k_v}) \neq \emptyset.
\end{equation}
The resulting connectivity is shown in Figure \ref{fig:hat_edge_definition}. By defining connectivity through this topological intersection, the pathfinding algorithm can transition fluidly between micro-scale navigation in complex terrain and macro-scale jumps across compressed regions, drastically reducing the total number of graph operations required to traverse the domain.
\subsection{Scale-Aware Edge Formulation}
Because the spatial distance between connected vertices is no longer constant, as it would be on uniform grids, the assignment of edge weights must be dynamically calculated based on the varying resolution levels. Let an edge connect vertex $u$ at resolution level $\ell_u$ to vertex $v$ at resolution level $\ell_v$. The physical distance $D(u,v)$ is computed using the Haversine formula based on the geographic coordinates of the respective cell centres. Haversine distance calculates the great-circle distance (``as the crow flies'') between two points on a sphere, accounting for the Earth's curvature. This is particularly important for large-scale routing, where planar approximations would introduce significant errors.
\par
Using the elevation values extracted from the multi-scale wavelet approximation, we evaluate the topographic slope $S(u,v)$ across this extended distance. The baseline traversal cost $C_{base}(u,v)$ is then strictly defined as the spatial distance divided by the scale-adjusted walking speed:
\begin{equation}
    C_{base}(u,v) = \frac{D(u,v)}{KS(S(u,v))}.
\end{equation}
This ensures that crossing a highly compressed, massive block is mathematically treated as a long-distance traversal rather than a single grid step.

\subsection{Topological Fidelity via Scale-Dependent Penalties}
When the wavelet transform consolidates a region into a coarse cell at level $\ell < L$, it filters out sub-threshold topographic variations. The algorithm models the path across this cell as a perfectly straight line on a smoothed surface. In reality, sub-scale roughness forces a traveller to deviate slightly, increasing the actual travel time. Without correction, the algorithm might exploit these mathematically smoothed blocks as artificial shortcuts.
\par
To enforce topological fidelity, we introduce a scale-dependent penalty factor $\alpha_\ell \geq 1.0$. The final modified edge cost $C(u,v)$ is calculated as:
\begin{equation}
C(u,v) = \alpha_{\ell_v} \cdot C_{base}(u,v)
\end{equation}
For fine-grained cells at maximum resolution $L$, we choose $\alpha_L = 1.0$ (no penalty). For larger, compressed blocks, $\alpha_\ell$ strictly increases, slightly inflating the calculated travel cost. This mathematical penalty dynamically compensates for the compressed sub-scale roughness, preventing the algorithm from unrealistically favouring coarse zones and ensuring high geometric accuracy across the multi-scale landscape.

\section{Case Studies and Results}
\label{sec:case_studies_and_results}

To validate the adaptive multi-scale routing framework, we benchmark its performance against the standard uniform Dijkstra algorithm using the 60 arc-second ETOPO dataset. The primary objective is to evaluate the trade-off between computational compression, measured as the reduction of spatial nodes, and the topological fidelity of the reconstructed pathways. 

For this validation, we selected two distinct geographic scenarios designed to test the algorithm to its limits: a macro-regional route across varied terrain (the Iberian Peninsula to the Western Alps) and a strictly micro-topographic challenge through dense mountain ranges (the Eastern Alps). The reference solutions were computed on a uniform grid containing $2^{20}$ vertices (prior to water mask application).

\subsection{Macro-Regional Routing: The Iberian Peninsula to the Western Alps}
The first case study models a macro-regional migration route spanning from the southern Iberian Peninsula ($38.6^\circ$N, $5.4^\circ$W) to the Western Alps ($45.0^\circ$N, $10.0^\circ$E), with a Haversine distance of approximately $1458$ km. This vast domain requires the algorithm to navigate a heterogeneous mix of flat coastal plains, the elevated Meseta Central, and significant mountainous barriers like the Pyrenees and western Alps.

We applied Best-$N$-term thresholding, which strictly retains only the $N$ largest wavelet detail coefficients, setting all others to zero. This allows us to precisely control the spatial compression rate. 

\paragraph{Graph Complexity vs.\ Routing Fidelity}
The results demonstrate the massive computational advantage of the adaptive framework. The compression rate acts as the primary indicator of efficiency. For example, using an $N_1$ wavelet basis with $N = 50,000$, the algorithm achieved a compression rate of $98.81\%$. Consequently, the total number of vertices ($|\mathcal{V}|$) the Dijkstra algorithm had to process dropped by over $80\%$, from approximately $285,000$ in the uncompressed reference grid to merely $52,899$.
\par
Operating an algorithm on an adaptive, multi-scale grid introduces additional computational overhead per node compared to a simple uniform grid. However, this initial investment pays off immensely at continental scales. By decisively shrinking the overall search space, the adaptive framework easily overcomes the processing limits of extremely large uniform grids.
\par
Crucially, despite discarding over $98\%$ of the mathematical detail coefficients, the physical distance of the resulting path remained highly stable, and the global routing topology was strictly preserved. As visualized in the resulting path models, the adaptive algorithm successfully navigated the massive, highly compressed spatial blocks in the plains of central Spain and France. However, upon reaching the Pyrenees, the underlying multi-resolution mesh dynamically forced the algorithm back to the highest geometric resolution. The adaptive paths identified and traversed the exact same major mountain corridors as the uncompressed reference solution.
\begin{figure}[htbp]
    \centering
    \begin{subfigure}[t]{0.45\textwidth}
        \includegraphics[width=\textwidth]{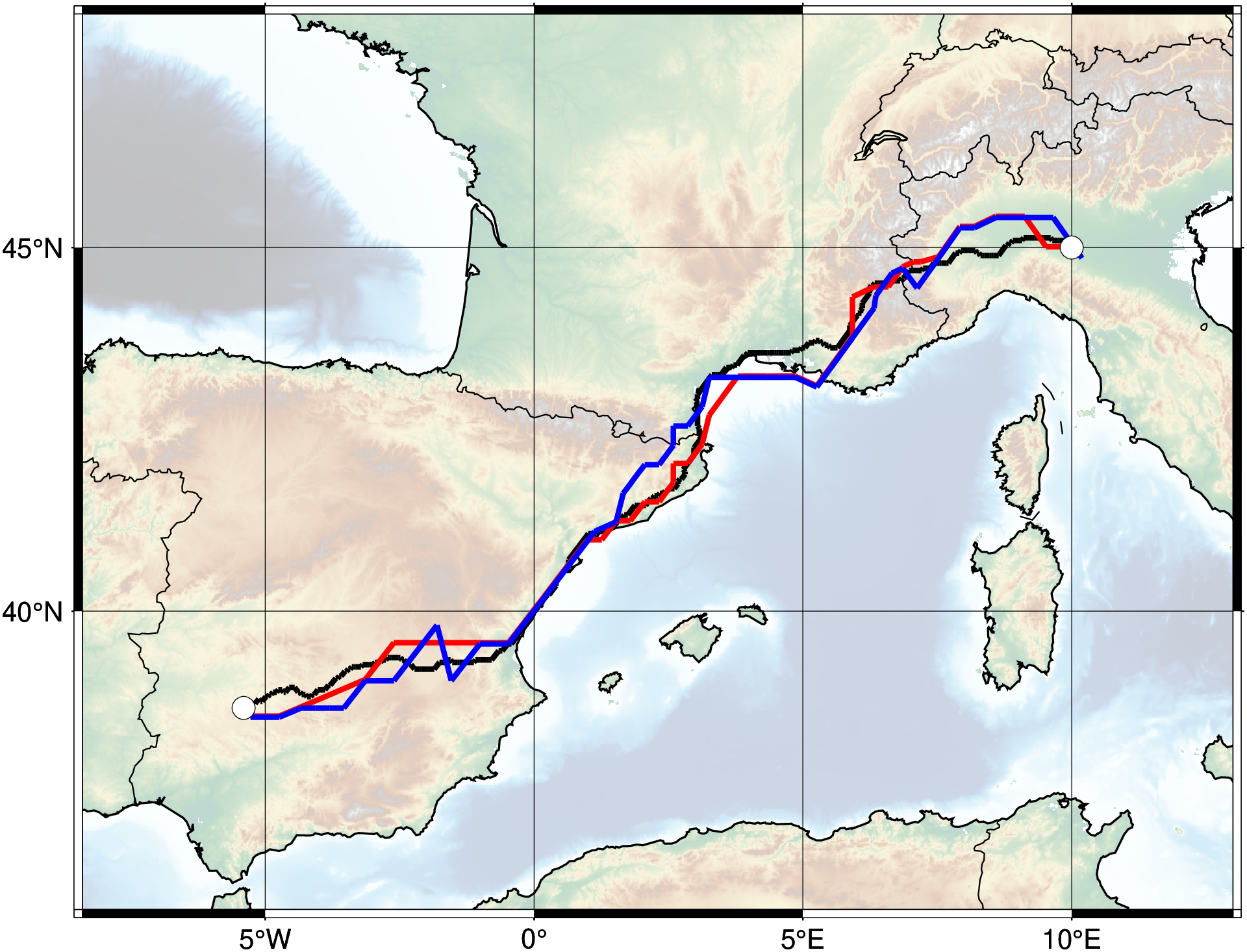}
        \caption{Paths computed using the adaptive $N_1$ (red) and $N_2$ (blue) wavelet bases at $N = 5000$ coefficients. The uncompressed reference path is shown in black. Even at high compression, the adaptive paths preserve the major mountain corridors, although some minor deviations are visible.}
        \label{fig:adaptive_paths_n5000}
    \end{subfigure}
    \hfill
    \begin{subfigure}[t]{0.45\textwidth}
        \includegraphics[width=\textwidth]{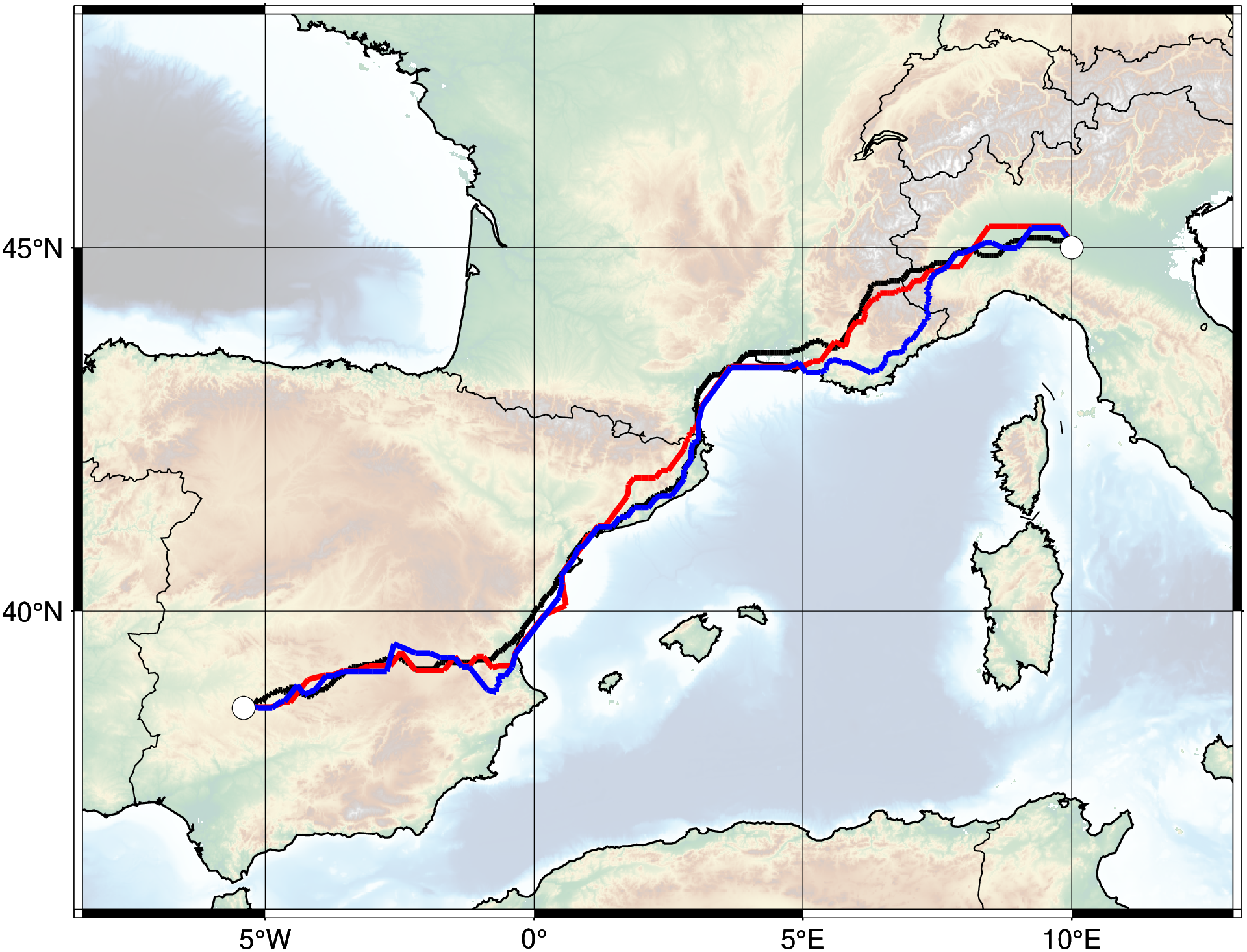}
        \caption{Paths computed using the adaptive $N_1$ (red) and $N_2$ (blue) wavelet bases at $N = 50,000$ coefficients. The uncompressed reference path is shown in black. At this moderate compression level, the adaptive paths closely align with the reference, demonstrating high topological fidelity.}
        \label{fig:adaptive_paths_n50000}
    \end{subfigure}
    \caption{Comparison of adaptive paths computed using different wavelet bases and compression levels. The uncompressed reference path, calculated at full resolution, is shown in black for context. (Figure adapted from \citep{brockmann_multigrid_2026})}
    \label{fig:adaptive_paths_comparison}
\end{figure}

\paragraph{Basis Selection: $N_1$ vs.\ $N_2$ Wavelets}
Comparing the multi-resolution bases revealed a clear performance dynamic. At extreme compression rates (e.g., $N = 5000$, $99.88\%$ compression), the piecewise-constant $N_1$ (Haar) basis proved more robust than the continuous piecewise-linear $N_2$ (Hat) basis. However, our benchmark showed that as more detail was retained (e.g., $N = 150,000$, $96.42\%$ compression), the $N_2$ basis leveraged its higher approximation order. The $N_2$ paths became nearly indistinguishable from the uncompressed reference, dropping the total cost error to $10.7\%$, compared to $18.6\%$ for the $N_1$ basis. Consequently, for heavily memory-constrained applications, the structurally simpler $N_1$ approach is highly efficient, whereas the $N_2$ basis is preferable when higher topological fidelity is required. How the compression rate and the resulting computational time are related is detailed in \cite{brockmann_2026_adaptive}. 

\subsection{Micro-Topographic Challenges: The Eastern Alps}
While macro-regional routing benefits immensely from spatial compression, highly rugged terrains push the mathematical limits of the adaptive mesh. To investigate these limits, the second case study evaluates a strictly localized route ($631.8$ km Haversine distance) through the dense, complex topography of the Eastern Alps. In such terrain, optimal movement is constrained to narrow valley floors.

\paragraph{The Valley Preservation Problem}
The reference path for this Alpine crossing exhibits a physical traversal distance of $1298.9$ km, more than double the straight-line Haversine distance. This highlights how the algorithm continuously zigzags through narrow valleys to avoid steep, energetically costly inclines. 
\par
\begin{figure}[htbp]
    \centering
    \begin{subfigure}[t]{0.3\textwidth}
        \includegraphics[width=\textwidth]{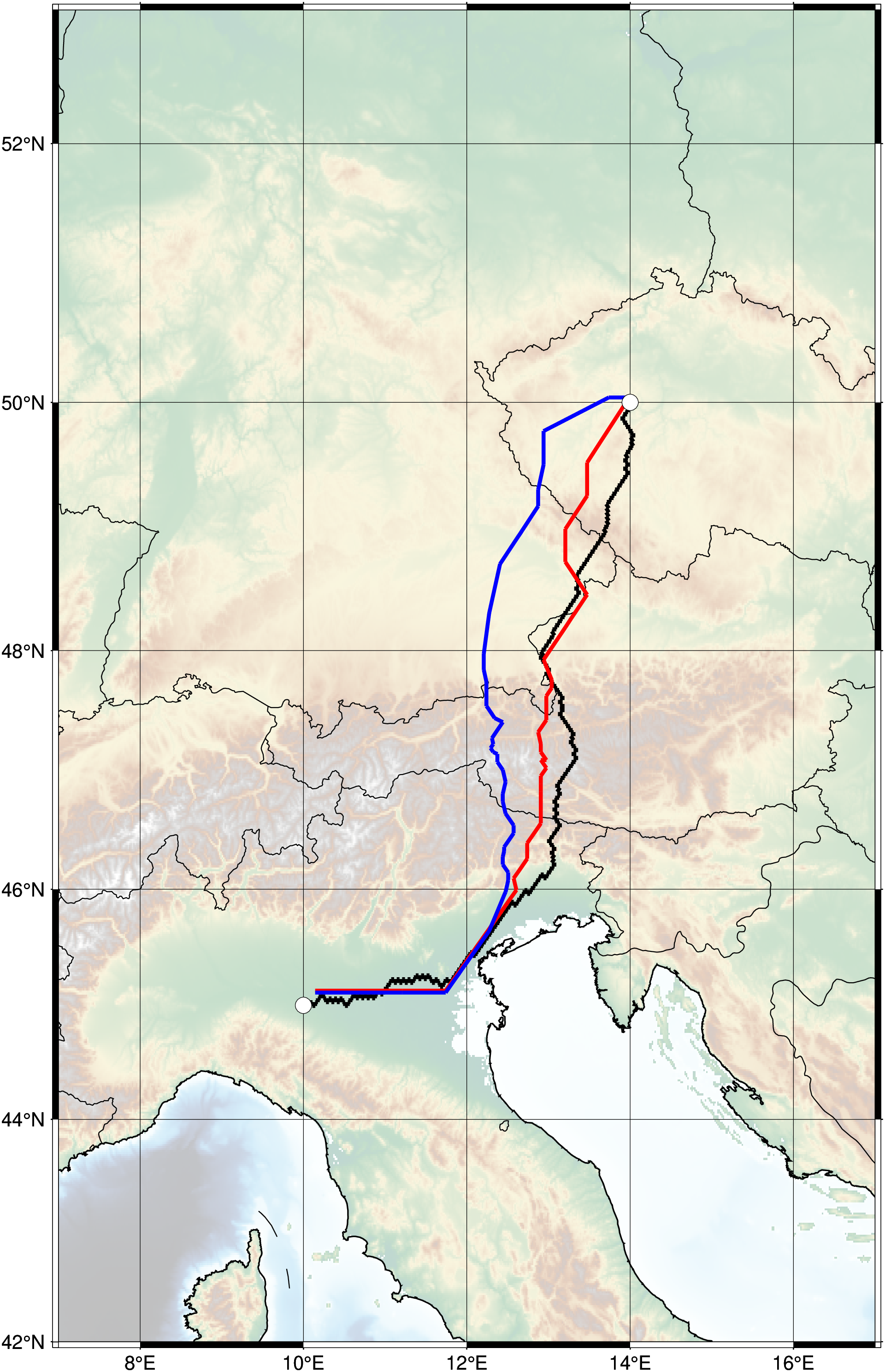}
        \caption{Paths computed using the adaptive $N_1$ (red) and $N_2$ (blue) wavelet bases at $N = 5000$ coefficients. The reference path on uncompressed data is shown in black. At this extreme compression level, the adaptive paths fail to preserve the narrow Alpine valleys, resulting in a path that fails to capture the mountainous terrain.}
        \label{fig:alpine_paths_n5000}
    \end{subfigure}\hfill
    \begin{subfigure}[t]{0.3\textwidth}
        \includegraphics[width=\textwidth]{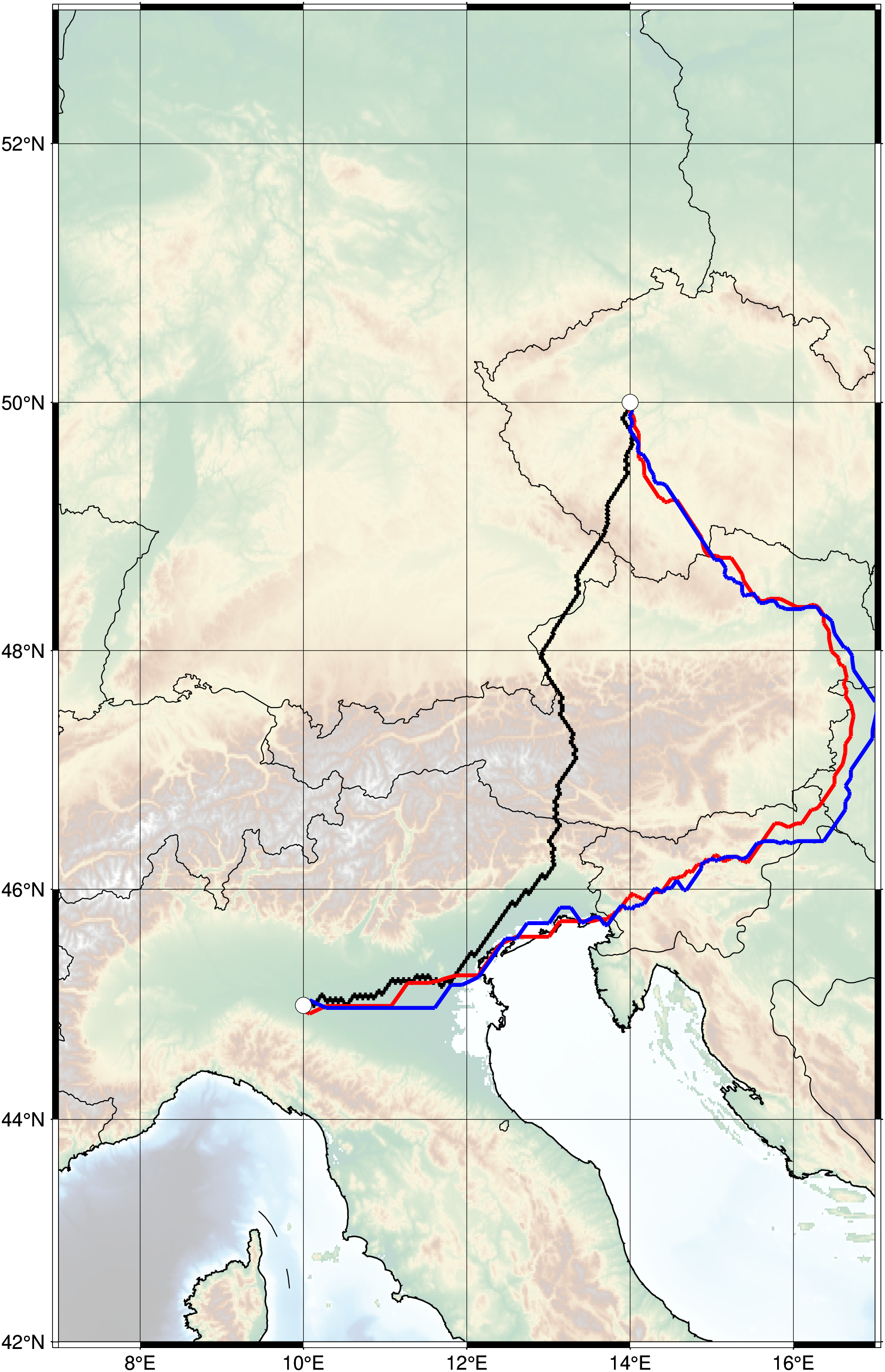}
        \caption{Paths computed using the adaptive $N_1$ (red) and $N_2$ (blue) wavelet bases at $N = 150,000$ coefficients. The reference path on uncompressed data is shown in black. At this moderate compression level, the adaptive paths recognize the mountains as a barrier but fail to preserve narrow passes, resulting in a detour around the mountain range.}
        \label{fig:alpine_paths_n150000}
    \end{subfigure}\hfill
    \begin{subfigure}[t]{0.3\textwidth}
        \includegraphics[width=\textwidth]{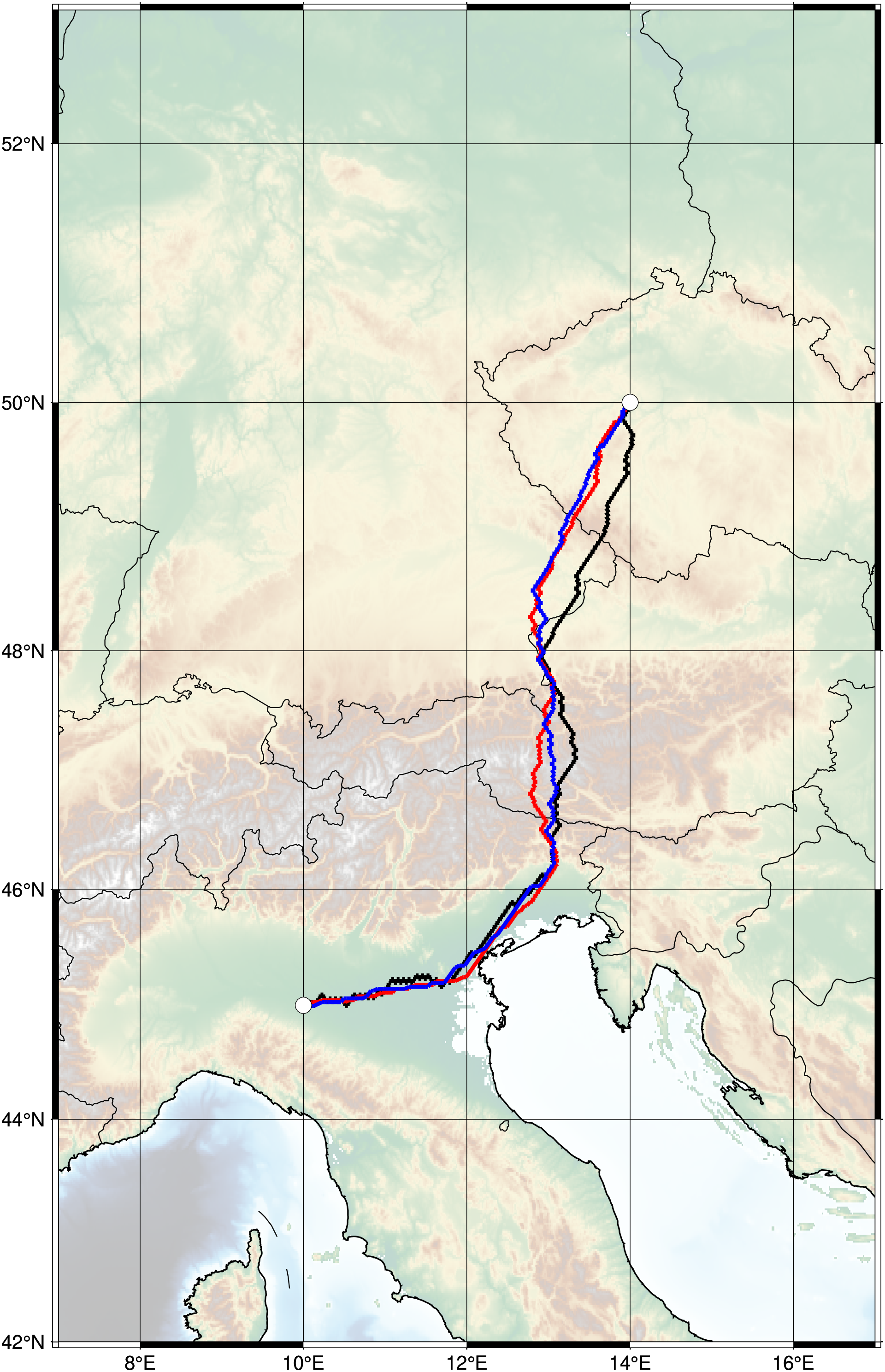}
        \caption{Paths computed using the adaptive $N_1$ (red) and $N_2$ (blue) wavelet bases at a compression rate of $\sim$ 32\%. The reference path on uncompressed data is shown in black. At this high detail level, the adaptive paths successfully preserve the narrow Alpine valleys, closely matching the reference path.}
        \label{fig:alpine_paths_r0_0005}
    \end{subfigure}
    \caption{Comparison of adaptive paths computed in the eastern Alps region, requiring the path to navigate through narrow valleys. The reference path, calculated at full resolution, is shown in black for context. (Figure adapted from \citep{brockmann_multigrid_2026})}
    \label{fig:alpine_paths_comparison}
\end{figure}

Applying adaptive compression in this environment reveals the structural vulnerability of barrier smearing. Wavelet compression fundamentally acts as a smoothing operator. When high compression rates are applied (e.g., $N = 5000$), the algorithm averages the high variations of the Alpine topography, flattening steep peaks and filling deep valleys. At these extreme thresholds, the physical traversal distance dropped dramatically to roughly $750$ km. The severe smoothing completely eradicated the Alpine barrier, causing the algorithm to perceive impassable mountains as a gentle incline and route an unrealistic, straight line across the topography. The resulting path is shown in Figure \ref{fig:alpine_paths_n5000}, where, despite following the correct general direction, the adaptive path visibly ignores the intricate valley corridors that define the true route.

\paragraph{Managing Topological Corridors}
Under moderate compression ($N = 150,000$), the mesh recovered enough detail to recognize the mountains as a formidable barrier. However, a secondary issue emerged: small geographic features, such as narrow mountain passes only a few pixels wide, were consolidated with adjacent steep slopes. With these vital topological corridors effectively ``closed'' by the compression smoothing, the routing algorithm was forced into massive geographic detours to circumvent the mountain range entirely, pushing the physical traversal distance back up to over $1100$ km, as shown in Figure \ref{fig:alpine_paths_n150000}.

This valley preservation problem was only resolved at the highest tested detail levels as shown in Figure \ref{fig:alpine_paths_r0_0005} (e.g., a compression rate of $\sim$ 32\%). At this threshold, the algorithm successfully preserved the narrow passes, allowing the adaptive paths to capture the intricate zigzagging through the Alpine valleys and converge with the uncompressed reference solution. This Alpine benchmark demonstrates that while adaptive meshes successfully reduce graph complexity in macro-regional scenarios, preserving the topological connectivity of highly constrained, sub-scale valleys remains a problem for which sufficient data resolution remains unavoidable.
\section{Discussion and Implications for Archaeological Spatial Modelling}
\label{sec:discussion_and_archaeological_implications}

The primary objective of this methodological study was to resolve the severe computational bottleneck that restricts large-scale, high-resolution mobility modelling in archaeological research. While we do not aim to reconstruct specific historical migration events in this paper, the proposed adaptive multi-scale routing framework provides a practical computational tool for future research within the HESCOR project. By significantly reducing algorithmic runtime, it allows researchers to readily incorporate high-resolution data into extensive spatial models without exceeding standard computing resources.

\subsection{Solving the Resolution-Scale Dilemma}
For spatial models that span a large geographic domain, the computational cost of evaluating high-resolution data fundamentally exceeds practical memory and runtime limits. Conversely, when downsampling is applied to reduce the computational load, the resulting loss of topographic detail can lead to fundamentally flawed path reconstructions.

Our benchmark across the Iberian Peninsula demonstrates that the adaptive wavelet framework effectively breaks this compromise. By reducing the processed vertex count by over $80\%$ while maintaining global route topology, the algorithm provides the computational efficiency required to evaluate massive spatial domains without discarding the underlying 60 arc-second ETOPO resolution. For projects analysing transcontinental dispersals processes or the spread of cultures, this framework enables the computation of APSP across a large set of archaeological sites, which would otherwise exceed practical memory and runtime limits on uniform grids.

\subsection{Topographic Fidelity and Human Movement}
The crucial advantage of the wavelet-based approach over uniform downsampling lies in its geometric preservation. Human movement is dictated by topographic extremes. Uniformly reducing grid resolution mathematically smooths these features out of existence, often leading routing algorithms to calculate paths over artificially flattened mountains rather than through the historically accurate, narrow valleys.
\par
The thresholding approach successfully mitigates this distortion in macro-regional models. By identifying areas of high local variance, the algorithm dynamically preserves the highest geometric resolution exactly where it is needed, at the structural choke points of the landscape. Consequently, we can trust that the macro-regional routes generated by this algorithm adhere to the actual topographic constraints of the environment, rather than artifacts of data compression.

\subsection{Strength of the Wavelet Approach}
Our benchmark scenarios across the Iberian Peninsula and the Eastern Alps demonstrate that the adaptive multi-scale routing framework effectively resolves the resolution-scale dilemma. Because wavelet thresholding allows for precise control over the retained topographic detail, it provides researchers with a flexible mechanism to balance computational efficiency against topological fidelity. For macro-regional models, the framework permits aggressive compression rates ($>95\%$) without sacrificing the global routing topology. Conversely, in micro-regional models constrained to rugged environments, a more conservative compression strategy ensures the preservation of narrow, sub-scale valleys that historically dictated human movement. By providing a mathematically optimized, adaptive representation of topography, this framework enables the efficient computation of extensive, high-resolution mobility networks that are otherwise computationally intractable.

\section{Conclusions}
\label{sec:conclusions}
In interdisciplinary network modelling within HESCOR, understanding raw material exchange or migration patterns often requires calculating the effective distance between hundreds of scattered settlement sites and resource origins. Furthermore, because human travel effort is inherently asymmetric, walking steep uphill gradients requires a different energy expenditure than walking downhill, the routing algorithm must independently evaluate paths in both directions for every pair of sites, effectively doubling the already substantial computational load.
\par
When calculating APSPs across continental domains, the use of high-resolution topographic data creates a severe computational bottleneck for standard Least-Cost Path algorithms. Conversely, uniformly downsampling these grids is methodologically flawed, as it inadvertently destroys critical topological corridors, such as narrow mountain passes, which historically dictated human movement.
\par
To resolve this resolution-scale dilemma, we introduced an adaptive multi-scale routing framework based on the Fast Wavelet Transform utilizing continuous piecewise-linear $N_2$ B-splines. Our approach dynamically compresses topographically homogeneous plains into large macroscopic blocks, while strictly preserving maximum geometric resolution in rugged, highly variable terrain.
\par
This adaptive spatial decomposition drastically reduces the routing graph size, thereby minimizing the computation time. Assisted by a scale-dependent penalty factor to prevent artificial routing shortcuts, the framework flawlessly maintains the topological fidelity of the pathways. Ultimately, this mathematically optimized tool overcomes the topographic bottleneck, making the generation of highly accurate, massive All-Pairs Shortest Paths matrices computationally feasible for future research. 

\section*{Acknowledgements}\vspace*{-1em}
This work was funded by the Ministry of Culture and Science of the State of North Rhine-Westphalia, 'HESCOR', PB22-081 ‘Profilbildung 2022’.
\nocite{*}
\printbibliography[heading=bibintoc, title={Bibliography}]

\end{document}